\documentclass[a4paper,11pt]{amsart}

\usepackage{cite}
\usepackage{url}

\usepackage{hyperref}

\usepackage{enumerate}

\newtheorem{theorem}{Theorem}[section]

\newtheorem{definition}[theorem]{Definition}

\numberwithin{equation}{section}

\newcommand{\Ocal} {{\mathcal O}}

\newcommand{\R}{\mathbb{R}}

\renewcommand{\P}{\mathbb{P}}
\newcommand{\E}{\mathbb{E}}

\newcommand{\sgn}{\operatorname{sgn}}
\renewcommand{\div}{\operatorname{div}}
\renewcommand{\epsilon}{\varepsilon}

\usepackage{enumerate}

\renewcommand{\div}{\operatorname{div}}

\usepackage{geometry}

\title[Stochastic evolution equations with nonlinear diffusivity]{Stochastic evolution equations with nonlinear diffusivity, recent progress and critical cases}

\author[I. Ciotir]{Ioana Ciotir}
\address{Normandie University\\
INSA de Rouen Normandie\\
LMI (EA 3226 - FR CNRS 3335)\\
76000 Rouen\\
France}
\email{ioana.ciotir@insa-rouen.fr}
\author[D. Goreac]{Dan Goreac}
\address{Universit\'e Laval\\
\'Ecole d'Actuariat\\
Facult\'e des sciences et de g\'enie\\
Pavillon Paul-Comtois\\
2425, rue de l'Agriculture\\
Local 4163\\
Qu\'ebec G1V 0A6\\
Canada\\
\& LAMA\\
Univ Gustave Eiffel, UPEM\\
Univ Paris Est Creteil, CNRS\\
77447 Marne-la-Vall\'ee\\
France}
\email{dan.goreac@act.ulaval.ca}
\author[J. M. T\"olle]{Jonas M. T\"olle}
\address{Aalto University\\
Department of Mathematics and Systems Analysis\\
PO Box 11100 (Otakaari 1, Espoo)\\
00076 Aalto\\
Finland}
\email{jonas.tolle@aalto.fi}

\date{\today}

\begin{document}

\begin{abstract}
This short survey article stems from recent progress on critical cases of stochastic evolution equations in variational formulation with additive, multiplicative or gradient noises. Typical examples appear as the limit cases of the stochastic porous medium equation, stochastic fast- and super fast-diffusion equations, self-organized criticality, stochastic singular $p$-Laplace equations, and the stochastic total variation flow, among others. We present several different notions of solutions, results on convergence of solutions depending on a parameter, and homogenization. Furthermore, we provide some references hinting at the recent progress in regularity results, long-time behavior, ergodicity, and numerical analysis.\end{abstract}



\maketitle

\section{Overview}

In this work,  we aim to provide a non-exhaustive survey on some recent results concerning, essentially, the existence and uniqueness of the solution for a certain class of stochastic partial differential equations (SPDEs) of monotone type with drift terms given by nonlinear second order quasi-linear diffusion operators. Our initial aim is to discuss critical cases related to the stochastic porous medium equation and the stochastic singular $p$-Laplace equation.
In this work, critical cases are meant to refer to equations at the boundaries of parameter families, limiting equations resulting from rescaling or homogenization, or equations with singular coefficients, or otherwise out of the range of the usual well-posedness theory.
For reasons of coherence, we also provide a non-complete number of references recollecting the recent progress in regularity results, long-time behavior, ergodicity, and numerical analysis for a wider class of nonlinear stochastic diffusion equations.

In abstract form, the equation can be written as
\[\partial_t u(t)=A(u(t))+B(u(t))\partial_t W(t),\quad u(0)=x,\]
where, 
\begin{itemize}
    \item $A$ is a nonlinear monotone second order differential operator, 
    \item $B$ is a zero or first order multiplicative noise coefficient, mainly assumed to be linear in this presentation,
    \item $W$ is cylindrical Gaussian noise, and 
    \item the initial datum $x$ is chosen in some separable Hilbert space.
\end{itemize}

\subsection{Stochastic porous medium type equation}

The stochastic theory of nonlinear partial differential equations has recently been intensively
studied for drift coefficients of the form $A(u)=\Delta \Psi(u) $, where $\Psi :\mathbb{\
R\rightarrow R}$ defined by $\Psi ( u) =u^{[m]}:=|u|^{m-1}u$ is a maximal
monotone operator. The equation is considered with additive, multiplicative, or divergence type noise, both for Gaussian space-time white noise and for colored noise.

\subsubsection{Existence results and regularity of the solution}

This type of equation has been studied using several approaches which are chosen depending on the regularity of the drift and the form of the noise.
We mention here the variational approach, the semi-group approach, and the martingale measure approach, but these are only a few methods available in the literature. 

The \textbf{variational approach} was introduced by Pardoux in his thesis (see \cite{pardoux}) and then developed by Krylov and Rozovski\u{\i} in \cite{krylov-rozowskii}. The main idea of the method lies in the construction of a Gelfand triple $V\subset H \subset V^{\ast}$ fitting the equation, as follows.
We consider a separable Hilbert space $H$, its dual $H^{\ast}$ and a Banach space $V$ such that $V \subset H$ is embedded continuously and densely. Then, for the dual $V^{\ast}$ of the latter, we have that $H^{\ast} \subset V^{\ast}$ continuously and densely, and, after identifying $H$ and $H^{\ast}$ via the Riesz isomorphism, we have that \[_{V^{\ast }}\left\langle z,v\right\rangle _{V}=\left\langle z,v\right\rangle_{H},\quad \forall z\in H,~v\in V.\] 
In the case of a porous medium equation defined on a bounded domain $\mathcal{O}$, we denote by 
$V:=\mathbb{L}^2(\mathcal{O})$ and by $H:=(H_{0}^{1}(\mathcal{O}))^{*}=H^{-1}(\mathcal{O})$. Since $H_{0}^{1}(\mathcal{O}) \subset \mathbb{L}^2(\mathcal{O})$ continuously and densely, we have also that 
$\mathbb{L}^2(\mathcal{O})\subset (H_{0}^{1}(\mathcal{O}))^{*}\textnormal{ continuously and densely.}$
Now we can identify $H$ with its dual $H^{*}$ via the corresponding Riesz isomorphism, so $H\equiv H^{*}$ in this sense and 
\begin{align}\label{Gelfand}
V=\mathbb{L}^2(\mathcal{O})\subset H\equiv H^{*} \subset (\mathbb{L}^2(\mathcal{O}))^{*}=V^{*}
\end{align}
continuously and densely. Note that $V^{*}$ is not the ordinary dual space of $V$ and cannot be identified with $\mathbb{L}^2(\mathcal{O})$ as a space of functions. (For more details, the reader is kindly referred to \cite{PR}).
We consider a stochastic differential equation on $H$ of the type%
\begin{equation}
dX\left( t\right) =A\left( X\left( t\right) \right) dt+B\left( X\left(
t\right) \right) dW\left( t\right) ,  \label{eq1}
\end{equation}%
where \begin{itemize}
    \item $W\left( t\right) ,$ $t\in \left[ 0,T\right] $ is a cylindrical Wiener
process defined on a complete probability space $\left( \Omega ,\mathcal{F},%
\mathbb{P}\right) $ endowed with a normal filtration $\left( \mathcal{F}_{t}\right)
_{t}$, and taking values in another separable Hilbert space $U$,
\item $B$ is
taking values in the space of Hilbert-Schmidt operators $L_{2}\left(
U;H\right) ,$ i.e. $B:V\rightarrow L_{2}\left( U;H\right) $, and\item  $%
A:V\rightarrow V^{\ast }$ is  a nonlinear operator satisfying suitable (monotonicity) conditions. 
\end{itemize}

In the variational approach, the concept of solution to equation (\ref{eq1})
is the following.

\begin{definition}[variational solution]
A solution to equation (\ref{eq1}) is a continuous $H-$valued $\left( 
\mathcal{F}_{t}\right) -$ adapted process $\left( X_{t}\right) _{t\in \left[
0,T\right] }$, such that its $dt\otimes \mathbb{P}$ equivalence class $%
\widehat{X}$ satisfies%
\begin{equation*}
\widehat{X}\in L^{2}\left( \left[ 0,T\right] \times \Omega ;V\right) \cap
L^{2}\left( \left[ 0,T\right] \times \Omega ;H\right)
\end{equation*}%
and $\mathbb{P-}$a.s.%
\begin{equation*}
X\left( t\right) =X\left( 0\right) +\int_{0}^{t}A\left( \overline{X}\left(
s\right) \right) ds+\int_{0}^{t}B\left( \overline{X}\left( s\right) \right)
dW\left( s\right) ,\quad t\in \left[ 0,T\right] ,
\end{equation*}%
where $\overline{X}$ is any $V-$valued progressively measurable $dt\otimes 
\mathbb{P-}$version of $\widehat{X}$.
\end{definition}

The reader is referred to Theorem 4.2.4 from \cite{PR} for an existence result of the
variational solution.

The \textbf{mild solution approach} inspired the analytically strong and analytically weak solution which have the following forms in the case of a porous medium equation 
\begin{equation}
\left\{ 
\begin{array}{l}
dX\left( t\right) = A(X\left( t\right))dt +B\left( X\left( t\right) \right) dW\left( t\right) ,\quad
t\in \left[ 0,T\right] , \\ 
X\left( 0\right) =x,%
\end{array}%
\right.  \label{eq2}
\end{equation}%
on a time interval $\left[ 0,T\right] ,$ where 
\begin{itemize}
    \item $A:D\left( A\right) \subset
H^{-1}(\mathcal{O})\rightarrow H^{-1}(\mathcal{O})$ is a porous medium operator, \item $x \in H^{-1}(\mathcal{O})$, \item $B:D\left( B\right) \subset H^{-1}(\mathcal{O})\rightarrow L_{2}\left(U;H^{-1}(\mathcal{O})\right) $ is a linear operator taking values in the space of
Hilbert-Schmidt operators $L_{2}\left( U;H^{-1}(\mathcal{O})\right)$, with $U$ being another separable Hilbert space, 
\item $W\left( t\right) ,$ $t\in \left[ 0,T\right] $ is a cylindrical Wiener
process defined on a complete probability space $\left( \Omega ,\mathcal{F},%
\mathbb{P}\right) $ with a normal filtration $\left( \mathcal{F}_{t}\right)
_{t}$, process taking its values in $U$.
\end{itemize}
\begin{definition}[analytically strong solution]
An $H-$valued predictable process $X\left( t\right) ,$ $t\in \left[ 0,T%
\right] $, which takes values in $D\left( A\right) \cap D\left( B\right) ,$ $%
\mathbb{P}-a.s.$ is called a strong solution to problem (\ref{eq2}) if, for
arbitrary $t\in \left[ 0,T\right]$, we have $\mathbb{P}-a.s.$ that 
\begin{equation*}
X\left( t\right) =x+\int_{0}^{t}A(X\left( s\right)) ds+\int_{0}^{t}B\left( X\left( s\right) \right)
dW\left( s\right) ,
\end{equation*}%
where the integrals on the right-hand side are well-defined.
\end{definition}
By extension, we have the following.
\begin{definition}[analytically weak solution]
An $H-$valued predictable process $X\left( t\right) $, $t\in \left[ 0,T%
\right] ,$ is said to be a weak solution to problem (\ref{eq2}) if $X$ takes its values in $D\left( B\right) $, and, for arbitrary $t\in \left[ 0,T\right]$, we have 
\begin{eqnarray*}
\left\langle X\left( t\right) ,\gamma \right\rangle &=&\left\langle x,\gamma
\right\rangle +\int_{0}^{t} \left\langle A(X\left( s\right)) ,\gamma \right\rangle  ds +\int_{0}^{t}\left\langle \gamma ,B\left( X\left( s\right) \right)
dW\left( s\right) \right\rangle ,\quad \mathbb{P}\text{-a.s.}
\end{eqnarray*}
for $\gamma$ being any eigenfunction of the homogeneous Dirichlet Laplace operator.
The integrals on the right-hand side are well-defined.
\end{definition}

The readers are referred to \cite{DPZ} for the classical existence results in this sense concerning general stochastic PDE and also for the relations between the different solutions. For the case of the porous medium equation see \cite{BDPR:16}. The $A$-related term can also be written in $L^2$, thus shifting the formulation to an explicit $\Psi$-term.

For the \textbf{martingale approach}, we shall consider the same framework as above,
and we define the martingale solution as follows.

\begin{definition}[martingale solution]
The system $\left( \Omega ,\mathcal{F},\left( \mathcal{F}_{t}\right) _{t},%
\mathbb{P},W,X\right) $ is said to be a martingale solution to (\ref{eq2})
if 
\begin{itemize}
    \item $\left( \Omega ,\mathcal{F},\left( \mathcal{F}_{t}\right) _{t},\mathbb{P}%
\right) $ is a filtered probability space on which $W$ is a cylindrical
Wiener process, and
\item $X$ is a $\left( \mathcal{F}_{t}\right) _{t}-$adapted
continuous process that satisfies $\mathbb{P}-a.s.$%
\begin{equation*}
X\left( t\right) =x+\int_{0}^{t}A\left( X\left( s\right) \right)
ds+\int_{0}^{t}B\left( X\left( s\right) \right) dW\left( s\right) .
\end{equation*}
\end{itemize}
\end{definition}

In the literature such solutions are also called \textbf{weak solution} (albeit the weakness is sought in a probabilistic sense). Since in PDE literature,
and in the present work the word weak has a different (PDE) meaning, in order to
avoid ambiguities, we use the term \textbf{martingale solution}. Note that, for
example, in the definition above, the solution is strong in the PDE sense.
When it is necessary, we shall emphasize the difference between the
\emph{analytically weak} solution (which is weak from the PDE point of view) and
the \emph{stochastically weak or martingale} solution (which is weak from the
stochastic point of view), we have the following.

\medskip

For subdifferential-type operators $A=-\partial\Phi$, where $\Phi:D(H)\subset H\to [0,+\infty)$ is a convex functional, we can define:

\begin{definition}[SVI solution]\label{defSVI}
The process $X$ is said to be a \textbf{stochastic variational inequality solution} (SVI-solution) to (\ref{eq2})
if $X$ is an $\left( \mathcal{F}_{t}\right)_{t}-$adapted
continuous process that satisfies
\[\E\int_0^T\Phi(X_t)\,dt<\infty,\]
and there exists a constant $C>0$ such that for every $t\in [0,T]$,
\begin{align*}
&\E\|X(t)-Z(t)\|_H^2+2\E\int_0^t\Phi(X(s))\,ds\\
\le& \|x-z\|^2_H+2\E\int_{0}^{t}\Phi(Z(s))\,ds\\
&+2\E\int_{0}^{t} \langle G(s),X(s)-Z(s)\rangle_H\,ds+C\E\int_{0}^{t} \| B(X(s))-H(s)\|^2_{L_2(U,H)}\,ds,
\end{align*}
for every sufficiently regular quadruple $(z,Z,G,H)$ that provides a strong solution to the stochastic equation in $H$
\[Z(t)=z+\int_0^t G(s)\,ds+\int_0^tH(s)\,dW(s),\quad t\in [0,T],\quad \P\text{-a.s.}\]
\end{definition}

The regularity of $(z,Z,G,H)$ has to be specified on a case-by-case basis depending on the equation considered. For the porous medium equation, for instance,
\[\Phi(u):=\frac{1}{m+1}\int_\Ocal |u(\xi)|^{m+1}\,d\xi.\]
See \cite{GT:16,GR:17,BR:13,GR:2015,GRW:24,T:20,CT:16}. This is the most robust notion of a solution under $\Gamma$- and Mosco convergence of the convex energy functionals $\Phi$, see \cite{GT:16}.

Keeping in mind that the porous medium operator can be written of the form $A(u)=\Delta \Psi(u) $, where $\Psi :\mathbb{\
R\rightarrow R}$ defined by $\Psi ( u) =u^{[m]}:=|u|^{m-1}u$, one can see that according to the value of $m$ the equation has several physical interpretations and specific technical difficulties.

\begin{enumerate}
    \item In the case $m>1,$ the corresponding equation describes the \textbf{slow diffusions}
(dynamics of fluids in porous media). The existence, uniqueness and
positivity, as well as finer behavior of the corresponding solution have already been studied in \cite{positivity,strong,BG2,strong2,RRW:2007} for the
stochastic case, as well as in the monograph \cite{BDPR:16}, and the references therein. For the deterministic case, see \cite{ARON} and \cite{VASQUEZ}. Well-posedness and improved regularity for general degenerate parabolic stochastic evolution equations has been studied in \cite{H:13,DMH:15,DHV:16}.
An alternative approach for linear multiplicative noise has been proposed in \cite{BR:15}.

\item The case $m\in ( 0,1) $ is relevant to the mathematical modeling
of \textbf{dynamics of an ideal gas in a porous medium}. Finite
time extinction is studied in 3 spatial dimensions for $m\in \left[ \frac{1}{5}
,1\right) $ in \cite{BDPR3}. See also \cite{BDPR12}.
Improvements of the Sobolev regularity for the stochastic fast diffusion equation with linear multiplicative noise appear in \cite{CGT:24}. Numerical analysis for porous medium/fast diffusion in a very weak formulation has been studied in \cite{ES:12}.

\item The stochastic diffusion equation for $m \in (-1,0)$ (corresponding to \textbf{super-fast diffusions}) with multiplicative It\^{o} noise was studied in \cite{eusuperfast}.
For the case $m\leq -1,$ it has been proved that, even in the deterministic
case, there exists no solution with finite mass, see \cite{VASQUEZ}.

\item The \textbf{logarithmic case} $\Psi (u) =\log(u)$ was studied for a multiplicative
noise in a bounded domain in \cite{LN,ln-noi}. 
\begin{enumerate}
    \item Note that, for positive solutions, the equation
can be seen as formally corresponding to the situation $m=0$,  in   \[ \operatorname{div}
(u^{[m-1]}\nabla u)=m^{-1} \Delta \left(  u^{[m]} \right) \] since \[\operatorname{div}
(u^{-1}\nabla u)=~\Delta  \log (u).\]

\item On the other hand, one can obtain the following interpretation for the limiting case $m=0$, if the limit \[u^{[m]}=|u|^{m-1}u\to|u|^{-1}u=:\sgn(u),\ m\searrow 0\] is interpreted point-wise. The resulting model is related to
the \textbf{singular case} $\Delta \sgn(u)$, where $\sgn(0)$ takes all values in $[-1,1]$, thus leading to a multi-valued differential inclusion, which have been discussed in
\cite{GT:14,RWX:24} for general differential operators replacing the Laplace operator. 

\item If interpreted from the point of view of $\Gamma$-convergence, one also obtains this limit, see \cite{GT:16},
This model is closely related to the \textbf{Bak-Tang-Wiesenfeld sand-pile model} of \textbf{self-organized
criticality}, see \cite{CCGS,criticality,BG1,N:21,B:13,BLGN:25}.
\end{enumerate}
\end{enumerate}
Further possibilities of interpretation also apply to the noise appearing in the systems.
\begin{enumerate}
\item The noise can be first interpreted in the sense of \textbf{It\^{o}
integral}. 
\begin{enumerate}
\item Usually,  we have the It\^{o} noise in a \textbf{multiplicative sense} which means
that $B$ is Lipschitz from $H$ into the space of Hilbert-Schmidt
operators in $H$.  The multiplicative noise acts proportionally with the solution, and, thus, can be interpreted as a \emph{relative error} contribution.  
\item Another option is the \textbf{additive noise} that does not depend
on the solution component $X$.  Concerning the interpretation, the additive noise can be seen as a stochastic source 
(or, otherwise, as an \emph{absolute error}). 
\item A more recent and original type of noise is the \textbf{gradient dependent noise}
which can be expressed as follows%
\begin{equation}
\left\{ 
\begin{array}{l}
dX\left( t\right) =A\left( X\left( t\right) \right) dt+B\left( \nabla
X\left( t\right) \right) dW\left( t\right) ,~t\in \left[ 0,T\right] , \\ 
X\left( 0\right) =x.%
\end{array}%
\right.  \label{Ee-grad}
\end{equation}%
The interest of such a noise is that it acts on each direction. 
From the point of view of the physical interpretation, the most commune example is the one of an operator $B$ of the form $B(x)=b\cdot \nabla x$ with $b$ a divergence free vector, and models the presence of a turbulence. The stochastic porous medium equation was studied with such a noise in the Itô case in \cite{euIto} proving the existence of a martingale, analytically weak solution. Finite time extinction for the 1D stochastic porous medium equation with transport noise has been proved in \cite{H:21}.

More recently, there has been a line of developments for gradient dependent noise in \cite{C:23,FGG:2022,DFG:2020,FG:24,FG:19}.

See \cite{CF:23,DKP:24} for recent developments for linear diffusion with gradient dependent noise.
\end{enumerate}

\item Another possibility is the use of \textbf{Stratonovich integral}, which provides an
alternative to It\^{o}'s approach.
As it is well known, the It\^{o} stochastic integral is convenient in some
problems because it provides a martingale. On the other hand, the Stratonovich
integral is particularly preferable in applications, since we can work with
such integrals in the same way as we would with the ordinary integrals of smooths
functions. One consequence is the fact that the Stratonovich integral is
\emph{stable with respect to changes in the random term}. 

If the stochastic integral is considered in the Stratonovich sense,  an equation of the type of
(\ref{eq2}) can be written as 
\begin{equation}
\left\{ 
\begin{array}{l}
dX\left( t\right) =A\left( X\left( t\right) \right) dt+B\left( X\left(
t\right) \right) \circ dW\left( t\right) ,~t\in \left[ 0,T\right], \\ 
X\left( 0\right) =x.%
\end{array}%
\right.  \label{EeS}
\end{equation}%
\end{enumerate}
All the different types of noise presented above (multiplicative, additive and divergence) can be considered in the sense of It\^{o} or
Stratonovich. The stochastic porous medium equation with divergence It\^{o} noise was studied in the particular case of a strongly monotone diffusion in \cite{Turra,euIto}.

- To the best of our knowledge, concerning the stochastic porous medium equation in \emph{unbounded domain}, the
only results known are \cite{BARBU20151024,RRW:2007,RWX:24,GRW:24,FG:2024,Francesco}.\\

- The case of \emph{logarithmic diffusion with a multiplicative Stratonovich noise} is treated in the very recent contribution \cite{ln-noi}.\\

- The stochastic porous medium equation with \emph{Robin boundary conditions, gravity-driven infiltration and multiplicative noise} appears in the preprint \cite{euDanRobin}.\\

- For improved regularity results for the stochastic porous medium and fast diffusion equations, see \cite{CGT:24,GR:17,BSW:22,G:12}.\\

- The porous medium equation with \emph{space-time white noise} has been studied in \cite{DGG:21}.

\subsubsection{Convergence results, homogenization, and ergodicity}

Let $\mathcal{O}$ be an open bounded domain of
$\mathbb{R}
^{d}$ $\left(  1\leq d\leq3\right)  $ with smooth boundary $\partial
\mathcal{O}$. We also consider the stochastic porous medium equation
\begin{equation}
\left\{
\begin{array}
[c]{ll}%
dX\left(  t\right)  -\Delta\Psi\left(  X\left(  t\right)  \right)  dt\ni
B\left(  X\left(  t\right)  \right)  dW\left(  t\right)  , &
\quad\text{in}\quad\left(  0,T\right)  \times\mathcal{O}\\
\Psi\left(  X\left(  t\right)  \right)  \ni0, & \quad\text{on}\quad\left(
0,T\right)  \times\partial\mathcal{O}\\
X\left(  0\right)  =x, & %
\end{array}
\right.  \label{ecu}%
\end{equation}
where \begin{itemize}
    \item $x$ is the initial data, and 
    \item $\Psi:\mathbb{R\rightarrow}2^{\mathbb{R}}$
is a maximal monotone (possibly multi-valued) graph with polynomial growth, and 
\item the noise is assumed to be multiplicative. 
\end{itemize}Such an equation satisfies the general existence result mentioned above and has an analytically weak solution as defined before.
We are interested in the continuous dependence of the solution as
function of $\Psi$ for the stochastic porous medium equation (\ref{ecu}). This problem is relevant in the \textbf{asymptotic analysis and approximation} of stochastic porous medium equations.

To this propose we consider a family of maximal monotone graphs $\left\{
\Psi^{\alpha}\right\}  _{\alpha>0},~\Psi$ and denote $A^{\alpha}=\Delta
\Psi^{\alpha}\left(  x\right),$  with polynomial growth independent of $\alpha$.

If we assume that $\Psi^{\alpha}\rightarrow\Psi$ as
$\alpha\rightarrow0$ in the graph sense, i. e.,
\[
\left(  1+\lambda\Psi^{\alpha}\right)  ^{-1}x\longrightarrow\left(
1+\lambda\Psi\right)  ^{-1}x,\quad\forall\lambda>0,\quad\forall x\in\mathbb{R}%
\]
for $\alpha\rightarrow0$, we have that the corresponding
solution $X^{\alpha}$ is convergent to the solution $X$ to (\ref{ecu}), i. e.,
\[
\underset{\alpha\rightarrow0}{\lim}\sup_{t\in [0,T]}\mathbb{E}\left\vert X^{\alpha}\left(
t\right)  -X\left(  t\right)  \right\vert _{H^{-1}\left(  \mathcal{O}\right)
}^{2}=0.
\] For more details see \cite{eutrotter1}.

This type of result has several applications as follows:

\begin{enumerate}
\item Let $\Psi^{\alpha}:\mathbb{R\rightarrow
}2^{\mathbb{R}}$ defined by%
\[
\Psi^{\alpha}\left(  x\right)  =\left\vert x\right\vert ^{\alpha}%
\operatorname{sign}x,\quad0\leq\alpha<1.
\]
This particular $\Psi^{\alpha}$ corresponds to the
\textbf{stochastic fast diffusion} equation and is relevant in \emph{plasma physics}.

The case $\alpha=0$ is relevant in stochastic models for \textbf{self-organized
criticality}, which describes the property of systems to have critical point as
\emph{attractor} and to \emph{reach spontaneously a critical state} (see \cite{criticality}).

As a consequence of the convergence result, we get the convergence of the solutions for $\alpha\rightarrow0$.
This proves that the solution to the equation
related to the model of self-organized criticality with stochastic
perturbation, can be approximated by a sequence of processes describing the
fast diffusion.

\item The diffusivity function $\Psi^{\alpha
}:\mathbb{R\rightarrow}2^{\mathbb{R}}$ of stochastic fast diffusion equation
can also be written as%
\[
\Psi^{\alpha}\left(  x\right)  =\left\vert x\right\vert ^{1-\alpha
}\operatorname{sign}x,\quad0<\alpha\leq1.
\]
In case when $\alpha$ is near $0,$ the corresponding equation can be regarded as a
\textbf{perturbation of stochastic heat equation}.

By the convergence result for $\alpha\rightarrow0$ we have that the solution $X^{\alpha}$ 
is convergent to the solution $X$ to the linear stochastic heat equations.

\item Let $\Psi:\mathbb{R\rightarrow}2^{\mathbb{R}}$
\ be a maximal monotone graph of the form%
\[
\Psi\left(  x\right)  =\left\{
\begin{array}
[c]{ll}%
\Psi_{1}\left(  x\right)  , & \quad\text{if }x<a\\
\left[  \Psi_{1}\left(  a\right)  ,\Psi_{2}\left(  a\right)  \right]  , &
\quad\text{if }x=a\\
\Psi_{2}\left(  x\right)  , & \quad\text{if }x>a
\end{array}
\right.  ,
\]
\ \ where $\Psi_{1}$ and $\Psi_{2}$ are continuous and monotonically
non-decreasing functions satisfying the assumption mentioned above.

We define the approximation
\[
\Psi^{\alpha}\left(  x\right)  =\left\{
\begin{array}
[c]{ll}%
\Psi_{1}\left(  x\right)  , & \text{if }x<a-\alpha\\
\Psi_{1}\left(  a-\alpha\right)  \dfrac{a+\alpha-x}{2\alpha}+\Psi_{2}\left(
a+\alpha\right)  \dfrac{a-\alpha-x}{-2\alpha}, & \text{if }a-\alpha\leq x\leq
a+\alpha\\
\Psi_{2}\left(  x\right)  , & \text{if }a+\alpha<x
\end{array}
\right.  .
\]

Using the convergence result, if we approximate a maximal monotone multivalued graph by a sequence of
continuous and monotonically non-decreasing functions, we have that the approximating solution is also convergent to the
solution to the limiting equation. This way some properties of the approximating solution
could be transferred to the solution to equation in the critical case.
\begin{enumerate}
\item As a particular case we have the \textbf{Heaviside step function}
\[
H\left(  x\right)  =\left\{
\begin{array}
[c]{ll}%
0, & \quad\text{if }x<0\\
\left[  0,1\right]  , & \quad\text{if }x=0\\
1, & \quad\text{if }x>0
\end{array}
\right.  ,
\]
which is relevant in the \textbf{anomalous (singular) diffusion} equation of the type
\[
dX\left(  t\right)  =\Delta\left(  H\left(  X\left(  t\right)  -x_{c}\right)
\right)  dt+B\left(  X\left(  t\right)  \right)  dW\left(  t\right)  ,
\]
with $x_{c}$ the critical value (see \cite{criticality}).

\item Another particular case is%
\[
\Psi\left(  x\right)  =\operatorname{sign}x=\left\{
\begin{array}
[c]{ll}%
\frac{x}{\left\vert x\right\vert },\medskip & \quad\text{if \ }x\neq0\\
\left[  -1,1\right]  , & \quad\text{if }x=0
\end{array}
\right.  ,
\]
which as mentioned above is relevant in the \textbf{stochastic models for
self-organized criticality}.
\end{enumerate}
\item Yet another application resides in \textbf{homogenization results} for systems driven by a (generic) maximal monotone function $\Phi$, and further involving a periodic contribution $a$. Define
\begin{equation*}
\Psi^{\alpha}\left( \xi,x\right) =a^{\alpha}\left( \xi\right) \Phi\left(
x\right) =a\left( \frac{\xi}{\alpha}\right) \Phi\left( x\right)
\end{equation*}
and $A^{\alpha}\left( X\right)
=-\Delta\Psi^{\alpha}\left( \cdot,X\left( \cdot\right) \right) $ with%
\begin{equation*}
D\left( A^{\alpha}\right) =\left\{ x\in H^{-1}\left( \mathcal{O}\right) \cap
L^{1}\left( \mathcal{O}\right) ;\Psi^{\alpha}\left( \cdot,X\left( \cdot\right)
\right) \in H_{0}^{1}\left( \mathcal{O}\right) \right\} .
\end{equation*}
The convergence result was generalized in \cite{eutrotter2} to the case of an operator $\Psi$ depending also on the space variable $\xi \in \mathcal{O}$ corresponding to a \emph{heterogeneous material}. This allows to obtain a homogenization result which means that the sequence of solutions $X^{\alpha}$ corresponding to the moisture level in heterogeneous media is converging to the solution to a homogeneous equation constructed with an operator 
\begin{equation*}
\Psi ^{\hom }\left( x\right) =M_{\mathcal{O}}\left( a\right) \Phi \left( x\right)
,\quad \forall x\in \mathbb{R}\text{.}
\end{equation*}%
where $M_{\mathcal{O}}\left( a\right) =\dfrac{1}{\left\vert \mathcal{O}\right\vert }%
\int_{\mathcal{O}}a\left( \xi\right) d\xi$ is the mean value of $a$ over $\mathcal{O}$, $\left\vert
\mathcal{O}\right\vert $\ is the Lebesgue measure of $\mathcal{O}$.\\

See also \cite{eutrotter3} for a homogenization result in the general variational framework and \cite{euDanTrotter} for a similar convergence result in the case of a Stefan-type equation with Robin boundary conditions.

\end{enumerate}

On the other hand, \textbf{ergodicity properties} for the porous medium and fast diffusion equations make the object of the contributions
\cite{DGT:20,GT:14,N:23}. 
\begin{enumerate}
    \item The long-time behavior of solutions to \textbf{stochastic porous media equations with nonlinear multiplicative noise on bounded domains} is studied in \cite{DGT:20}, establishing the existence and uniqueness of invariant measures.
    \item Abstract variational existence and uniqueness results are established for \textbf{multi-valued, monotone, non-coercive stochastic evolution inclusions} in Hilbert spaces in \cite{GT:14}.
    \item The long-time behavior of solutions to \textbf{generalized stochastic porous media equations on bounded domains with Dirichlet boundary conditions} is analyzed, establishing existence and uniqueness of invariant measures in \cite{N:23}.
    \end{enumerate}

\subsection{Stochastic singular $p$-Laplace-type equation}

The \textbf{stochastic $p$-Laplace equation} has the drift contribution involving
\[A(u)=\div(|\nabla u|^{p-2}\nabla u),\quad p>1,\]
see \cite{DB:93} for the deterministic case. For simplicity, impose Dirichlet boundary conditions on a bounded domain $\Ocal\subset\R^d$ with Lipschitz boundary, or, to reduce boundary issues completely, periodic boundary conditions on a cube.

 The stochastic $p$-Laplace has been first studied in \cite{L:09}, and provides a standard example of a monotone drift SPDE, see e.g. \cite{PR}. For linear multiplicative noise, it has been studied via the alternative approach proposed in \cite{BR:15}.
 The singular case $p\in (1,2)$ lacks good coercivity estimates in higher spatial dimensions due to restriction $p\ge\frac{2d}{d+2}$ of the Sobolev embedding $W^{1,p}\subset L^2$.
 
We distinguish six cases with particular relevance to the choice of the spaces on which the arguments are developed.

\begin{enumerate}
\item \textbf{Case $p>2$}
 
 This case is called the \textbf{degenerate case}. The Gelfand triple is \[V=W^{1,p}_0(\Ocal),\ H=L^2(\Ocal),\ V^\ast=W^{-1,q}(\Ocal),\] where $\frac{1}{p}+\frac{1}{q}=1$. By standard embeddings, $V\subset L^p(\Ocal)\subset H$, densely and continuously. See \cite{DHW:23,BHL:21} for numerical analysis of the stochastic $p$-Laplace system for $p>2$. Improved temporal regularity in this case has been obtained in \cite{W:23}. See \cite{SZ:23} for existence and uniqueness on the whole space $\R^d$. The case of stochastic $p(t,x)$-Laplace with variable exponents has been initiated in \cite{BVWZ:13}. $L^1$-initial data have been treated in \cite{SZ:22,SZ:21}.

\item \textbf{Case $p=2$}

This case reduces to the \textbf{linear drift stochastic heat equation}, $A=\Delta$. The Gelfand triple is \[V=H^1_0(\Ocal)=W^{1,2}_0(\Ocal),\ H=L^2(\Ocal),\ V^\ast=H^{-1}(\Ocal)=(H^1_0(\Ocal))^\ast.\] This case has been well-studied, see \cite{DPZ} and the references therein.\\
As the equation has a linear principal term, \emph{fixed point} methods apply.\\
One can also weaken the assumptions on the noise coefficient $B$ and obtain \textbf{mild solutions}.

\item \textbf{Case $1\vee\frac{2d}{d+2}\le p<2$}

This case is called the \textbf{subcritical singular case}. The Gelfand triple is \[V=W^{1,p}_0(\Ocal),\ H=L^2(\Ocal),\ V^\ast=W^{-1,q}(\Ocal),\] where $\frac{1}{p}+\frac{1}{q}=1$.\\
By the Sobolev embedding theorem, $V\subset H$ densely and continuously. If $p> 1\vee\frac{2d}{d+2}$, the embedding is compact by the Rellich-Kondrachov theorem.\\
Here, the \emph{standard theory} of monotone variational SPDEs works, see \cite{PR}.

\item \textbf{Case $1<p<\frac{2d}{d+2}$}

This case is called the \textbf{critical singular case}. In this case, it is more useful to consider the Gelfand triple \[V=H^1_0(\Ocal)=W^{1,2}_0(\Ocal),\ H=L^2(\Ocal),\ V^\ast=H^{-1}(\Ocal)=(H^1_0(\Ocal))^\ast,\] together with \textbf{generalized solutions}, that is, \emph{solutions obtained by limits of approximations} (SOLA), or SVI solutions, see \cite{GT:14,GR:17,GT:16,T:20}, and refer to Definition \ref{defSVI}.

\item\textbf{$1\vee (2-\frac{4}{d})<p<2$}
 
 In this case, we consider the same Gelfand triple as in the previous case \[V=H^1_0(\Ocal)=W^{1,2}_0(\Ocal),\ H=L^2(\Ocal),\ V^\ast=H^{-1}(\Ocal)=(H^1_0(\Ocal))^\ast,\] together with \textbf{generalized solutions} or \textbf{SVI solutions}.\\
 In this case, one can obtain \emph{second order estimates} and \emph{improved regularity} for the solutions and the decay behavior of the dynamics, see \cite{SST:23}.

\item\textbf{Case $p=1$}

In this case the term
\[A(u)=\div\left(\frac{\nabla u}{|\nabla u|}\right)\]
is not anymore pointwise well-defined.\\
Firstly, the divergence has to be interpreted in the \textbf{distributional sense}. Second, the singular diffusivity has to be extended to a multi-valued operator, as it cannot anymore be ensured that the diffusivity is a.e. zero, when $\nabla u=0$ a.e. Fortunately, by \emph{convex analysis and variational methods}, one can make sense of $A$ as the subdifferential of the relaxation of the total variation of the gradient measure of $H^1$-functions in the space $BV$. Then, $A$ is called the $1$-Laplace and the equation is called \textbf{total variation flow}. As $BV$ is not reflexive and is (compactly) embedded to $L^2$ only in dimensions $d=1,2$, (only in dimension $d=2$) one cannot use Gelfand triples anymore, and has to rely on \textbf{limit solutions} \cite{GT:14} or the more robust approach of \textbf{SVI solutions}, see \cite{BR:13,GR:17,GT:16}.
SVI solutions for the stochastic total variation flow have been extended to linear gradient Stratonovich noise in \cite{CT:16,T:20,T:18}.
In \cite{BRW1}, the authors propose a convergent numerical scheme for the stochastic total variation flow.
\end{enumerate}

\textbf{Further results}

\textbf{Ergodicity, stability and moment estimates} for the stochastic singular $p$-Laplace equation have been studied in \cite{GT:14,GT:16,LT:11,SST:23,GT-siam}.
\begin{enumerate}
    \item Ergodicity is established in \cite{GT:16} for local and nonlocal stochastic singular $p$-Laplace equations across all spatial dimensions, including multivalued cases, with convergence of invariant measures demonstrated under appropriate conditions.
    \item Existence and uniqueness of invariant measures are established in \cite{LT:11} for stochastic evolution equations with weakly dissipative drifts using new decay estimates and applications to singular $p$-Laplace and fast diffusion equations.
    \item \cite{SST:23} presents stability, long-time behavior, and moment estimates for stochastic singular $\Phi$-Laplace equations with additive Wiener noise and singular drift.
    \item In \cite{GT-siam}, ergodicity and convergence of invariant measures are proven for stochastic local and nonlocal singular $p$-Laplace equations, including multivalued cases, without spatial dimension restrictions and for all $p\in[1,2)$.
\end{enumerate}


\begin{thebibliography}{10}

\bibitem{ARON}
D.~G. Aronson.
\newblock {\em The porous medium equation}, pages 1--46.
\newblock Springer Berlin Heidelberg, Berlin, Heidelberg, 1986.

\bibitem{BLGN:25}
L.~Ba{\v{n}}as, B.~Gess, and M.~Neu{\ss}.
\newblock Stochastic partial differential equations arising in self-organized
  criticality.
\newblock {\em Ann. Appl. Probab.}, 35(1):481--522, 2025.

\bibitem{LN}
V.~Barbu.
\newblock The fast logarithmic equation with multiplicative {G}aussian noise.
\newblock {\em Ann. Univ. Buchar. Math. Ser.}, 2, 01 2012.

\bibitem{B:13}
V.~Barbu.
\newblock Self-organized criticality of cellular automata model; absorbtion in
  finite-time of supercritical region into the critical one.
\newblock {\em Math. Methods Appl. Sci.}, 36(13):1726--1733, 2013.

\bibitem{BDPR3}
V.~Barbu, G.~Da~Prato, and M.~R\"{o}ckner.
\newblock Finite time extinction for solutions to fast diffusion stochastic
  porous media equations.
\newblock {\em C. R. Math. Acad. Sci. Paris}, 347(1-2):81--84, 2009.

\bibitem{BDPR:16}
V.~Barbu, G.~Da~Prato, and M.~R\"{o}ckner.
\newblock {\em Stochastic porous media equations}, volume 2163 of {\em Lecture
  Notes in Mathematics}.
\newblock Springer, [Cham], 2016.

\bibitem{positivity}
V.~Barbu, G.~Da~Prato, and M.~Röckner.
\newblock Existence and uniqueness of nonnegative solutions to the stochastic
  porous media equation.
\newblock {\em Indiana Univ. Math. J.}, 57, 04 2007.

\bibitem{BDPR12}
V.~Barbu, G.~{Da Prato}, and M.~Röckner.
\newblock Finite time extinction of solutions to fast diffusion equations
  driven by linear multiplicative noise.
\newblock {\em J. Math. Anal. Appl.}, 389(1):147--164, 2012.

\bibitem{strong}
V.~Barbu, G.~D. Prato, and M.~Röckner.
\newblock {Existence of strong solutions for stochastic porous media equation
  under general monotonicity conditions}.
\newblock {\em Ann. Probab.}, 37(2):428 -- 452, 2009.

\bibitem{BR:13}
V.~Barbu and M.~R\"{o}ckner.
\newblock Stochastic variational inequalities and applications to the total
  variation flow perturbed by linear multiplicative noise.
\newblock {\em Arch. Ration. Mech. Anal.}, 209(3):797--834, 2013.

\bibitem{BR:15}
V.~Barbu and M.~R\"{o}ckner.
\newblock An operatorial approach to stochastic partial differential equations
  driven by linear multiplicative noise.
\newblock {\em J. Eur. Math. Soc. (JEMS)}, 17(7):1789--1815, 2015.

\bibitem{criticality}
V.~Barbu and M.~Röckner.
\newblock Stochastic porous media equations and self-organized criticality:
  Convergence to the critical state in all dimensions.
\newblock {\em Comm. Math. Phys.}, 311, 02 2011.

\bibitem{BARBU20151024}
V.~Barbu, M.~Röckner, and F.~Russo.
\newblock Stochastic porous media equations in $\mathbb{R}^d$.
\newblock {\em J. Math. Pures Appl.}, 103(4):1024--1052, 2015.

\bibitem{BVWZ:13}
C.~Bauzet, G.~Vallet, P.~Wittbold, and A.~Zimmermann.
\newblock On a {$p(t,x)$}-{L}aplace evolution equation with a stochastic force.
\newblock {\em Stoch. Partial Differ. Equ. Anal. Comput.}, 1(3):552--570, 2013.

\bibitem{BRW1}
L.~Ba\v{n}as, M.~R\"ockner, and A.~Wilke.
\newblock Convergent numerical approximation of the stochastic total variation
  flow.
\newblock {\em Stoch. Partial Differ. Equ. Anal. Comput.}, 9(2):437--471, 2021.

\bibitem{BHL:21}
D.~Breit, M.~Hofmanov\'{a}, and S.~Loisel.
\newblock Space-time approximation of stochastic {$p$}-{L}aplace-type systems.
\newblock {\em SIAM J. Numer. Anal.}, 59(4):2218--2236, 2021.

\bibitem{BSW:22}
S.~Bruno, B.~Gess, and H.~Weber.
\newblock Optimal regularity in time and space for stochastic porous medium
  equations.
\newblock {\em Ann. Probab.}, 50(6):2288--2343, 2022.

\bibitem{CCGS}
J.~Carlson, J.~Chayes, E.~Grannan, and G.~Swindle.
\newblock Self-organized criticality in sandpiles: nature of the critical
  phenomenon.
\newblock {\em Phys. Rev. A}, 42(4):2467--2470, 1990.

\bibitem{eutrotter1}
I.~Ciotir.
\newblock A {T}rotter type result for the stochastic porous media equations.
\newblock {\em Nonlinear Anal.}, 71:5606--5615, 2009.

\bibitem{eutrotter2}
I.~Ciotir.
\newblock Convergence of solutions for the stochastic porous media equations
  and homogenization.
\newblock {\em J. Evol. Equ.}, 11:339--370, 2011.

\bibitem{eutrotter3}
I.~Ciotir.
\newblock A {T}rotter-type theorem for nonlinear stochastic equations in
  variational formulation and homogenization.
\newblock {\em Differ. Integral Equ.}, 24(3--4):371--388, 2011.

\bibitem{eusuperfast}
I.~Ciotir.
\newblock {Existence and uniqueness of the solution for stochastic super-fast
  diffusion equations with multiplicative noise}.
\newblock {\em {Aust. J. Math. Anal. Appl.}}, 452(1):595--610, Aug. 2017.

\bibitem{euIto}
I.~Ciotir.
\newblock Stochastic porous media equations with divergence {I}t\^{o} noise.
\newblock {\em Evol. Equ. Control Theory}, 9(2):375--398, 2020.

\bibitem{ln-noi}
I.~Ciotir, R.~Fukuizumi, and D.~Goreac.
\newblock The stochastic fast logarithmic equation in {$\mathbb{R}^d$} with
  multiplicative {S}tratonovich noise.
\newblock {\em J. Math. Anal. Appl.}, 542(1):Paper No. 128786, 25, 2025.

\bibitem{euDanRobin}
I.~Ciotir, D.~Goreac, J.~Li, and Y.~Peng.
\newblock Trotter-type results for a stochastic {S}tefan-type system with
  {R}obin boundary conditions.
\newblock {\em to appear in Comput. Appl. Math.}, 2025+.

\bibitem{euDanTrotter}
I.~Ciotir, D.~Goreac, J.~Li, and A.~Tonnoir.
\newblock Stochastic porous media equation with robin boundary conditions,
  gravity-driven infiltration and multiplicative noise.
\newblock {\em J. Differential Equations}, 438:113363, 2025.

\bibitem{CGT:24}
I.~Ciotir, D.~Goreac, and J.~M. T\"{o}lle.
\newblock Improved regularity for the stochastic fast diffusion equation.
\newblock {\em Electron. Commun. Probab.}, 29:Paper No. 5, 7, 2024.

\bibitem{CT:16}
I.~Ciotir and J.~M. T\"{o}lle.
\newblock Nonlinear stochastic partial differential equations with singular
  diffusivity and gradient {S}tratonovich noise.
\newblock {\em J. Funct. Anal.}, 271(7):1764--1792, 2016.

\bibitem{C:23}
A.~Clini.
\newblock Porous media equations with nonlinear gradient noise and {D}irichlet
  boundary conditions.
\newblock {\em Stochastic Process. Appl.}, 159:428--498, 2023.

\bibitem{CF:23}
F.~Cornalba and J.~Fischer.
\newblock The {D}ean-{K}awasaki equation and the structure of density
  fluctuations in systems of diffusing particles.
\newblock {\em Arch. Ration. Mech. Anal.}, 247(5):Paper No. 76, 59, 2023.

\bibitem{DPZ}
G.~Da~Prato and J.~Zabczyk.
\newblock {\em Stochastic equations in infinite dimensions}, volume 152 of {\em
  Encyclopedia of Mathematics and its Applications}.
\newblock Cambridge University Press, Cambridge, second edition, 2014.

\bibitem{DGG:21}
K.~Dareiotis, M.~Gerencs\'er, and B.~Gess.
\newblock Porous media equations with multiplicative space-time white noise.
\newblock {\em Ann. Inst. Henri Poincar\'e{} Probab. Stat.}, 57(4):2354--2371,
  2021.

\bibitem{DGT:20}
K.~Dareiotis, B.~Gess, and P.~Tsatsoulis.
\newblock Ergodicity for stochastic porous media equations with multiplicative
  noise.
\newblock {\em SIAM J. Math. Anal.}, 52(5):4524--4564, 2020.

\bibitem{DMH:15}
A.~Debussche, S.~de~Moor, and M.~Hofmanov\'{a}.
\newblock A regularity result for quasilinear stochastic partial differential
  equations of parabolic type.
\newblock {\em SIAM J. Math. Anal.}, 47(2):1590--1614, 2015.

\bibitem{DHV:16}
A.~Debussche, M.~Hofmanov\'{a}, and J.~Vovelle.
\newblock Degenerate parabolic stochastic partial differential equations:
  quasilinear case.
\newblock {\em Ann. Probab.}, 44(3):1916--1955, 2016.

\bibitem{DB:93}
E.~DiBenedetto.
\newblock {\em Degenerate parabolic equations}.
\newblock Universitext. Springer-Verlag, New York, 1993.

\bibitem{DHW:23}
L.~Diening, M.~Hofmanov\'{a}, and J.~Wichmann.
\newblock An averaged space-time discretization of the stochastic
  {$p$}-{L}aplace system.
\newblock {\em Numer. Math.}, 153(2-3):557--609, 2023.

\bibitem{DFG:2020}
N.~Dirr, B.~Fehrman, and B.~Gess.
\newblock Conservative stochastic {PDE} and fluctuations of the symmetric
  simple exclusion process.
\newblock {\em Preprint}, pages 1--39, 2020.
\newblock https://arxiv.org/abs/2012.02126.

\bibitem{DKP:24}
A.~Djurdjevac, H.~Kremp, and N.~Perkowski.
\newblock Weak error analysis for a nonlinear {SPDE} approximation of the
  {D}ean-{K}awasaki equation.
\newblock {\em Stoch. Partial Differ. Equ. Anal. Comput.}, 12(4):2330--2355,
  2024.

\bibitem{ES:12}
E.~Emmrich and D.~{\v{S}}i{\v{s}}ka.
\newblock Full discretization of the porous medium/fast diffusion equation
  based on its very weak formulation.
\newblock {\em Commun. Math. Sci.}, 10(4):1055--1080, 2012.

\bibitem{FG:19}
B.~Fehrman and B.~Gess.
\newblock Well-posedness of nonlinear diffusion equations with nonlinear,
  conservative noise.
\newblock {\em Arch. Ration. Mech. Anal.}, 233(1):249--322, 2019.

\bibitem{FG:2024}
B.~Fehrman and B.~Gess.
\newblock Conservative stochastic {PDE}s on the whole space.
\newblock {\em Preprint}, pages 1--26, 2024.
\newblock https://arxiv.org/abs/2410.00254.

\bibitem{FG:24}
B.~Fehrman and B.~Gess.
\newblock Well-posedness of the {D}ean-{K}awasaki and the nonlinear
  {D}awson-{W}atanabe equation with correlated noise.
\newblock {\em Arch. Ration. Mech. Anal.}, 248(2):Paper No. 20, 60, 2024.

\bibitem{FGG:2022}
B.~Fehrman, B.~Gess, and R.~S. Gvalani.
\newblock Ergodicity and random dynamical systems for conservative {SPDE}s.
\newblock {\em Preprint}, pages 1--66, 2024.
\newblock https://arxiv.org/abs/2206.14789.

\bibitem{BG2}
B.~Gess.
\newblock Finite speed of propagation for stochastic porous media equations.
\newblock {\em SIAM J. Math. Anal.}, 45, 10 2012.

\bibitem{G:12}
B.~Gess.
\newblock Strong solutions for stochastic partial differential equations of
  gradient type.
\newblock {\em J. Funct. Anal.}, 263(8):2355--2383, 2012.

\bibitem{BG1}
B.~Gess.
\newblock Finite time extinction for stochastic sign fast diffusion and
  self-organized criticality.
\newblock {\em Commun. Math. Phys.}, 335, 10 2013.

\bibitem{GR:2015}
B.~Gess and M.~R\"{o}ckner.
\newblock Singular-degenerate multivalued stochastic fast diffusion equations.
\newblock {\em SIAM J. Math. Anal.}, 47(5):4058--4090, 2015.

\bibitem{GR:17}
B.~Gess and M.~R\"{o}ckner.
\newblock Stochastic variational inequalities and regularity for degenerate
  stochastic partial differential equations.
\newblock {\em Trans. Amer. Math. Soc.}, 369(5):3017--3045, 2017.

\bibitem{GRW:24}
B.~Gess, M.~R\"{o}ckner, and W.~Wu.
\newblock S{VI} solutions to stochastic nonlinear diffusion equations on
  general measure spaces.
\newblock {\em J. Evol. Equ.}, 24(4):Paper No. 94, 37, 2024.

\bibitem{GT:14}
B.~Gess and J.~M. T\"{o}lle.
\newblock Multi-valued, singular stochastic evolution inclusions.
\newblock {\em J. Math. Pures Appl. (9)}, 101(6):789--827, 2014.

\bibitem{GT-siam}
B.~Gess and J.~M. T\"{o}lle.
\newblock Ergodicity and local limits for stochastic local and nonlocal
  {$p$}-{L}aplace equations.
\newblock {\em SIAM J. Math. Anal.}, 48(6):4094--4125, 2016.

\bibitem{GT:16}
B.~Gess and J.~M. T\"{o}lle.
\newblock Stability of solutions to stochastic partial differential equations.
\newblock {\em J. Differential Equations}, 260(6):4973--5025, 2016.

\bibitem{H:21}
S.~Hensel.
\newblock Finite time extinction for the 1{D} stochastic porous medium equation
  with transport noise.
\newblock {\em Stoch. Partial Differ. Equ. Anal. Comput.}, 9(4):892--939, 2021.

\bibitem{H:13}
M.~Hofmanov\'{a}.
\newblock Degenerate parabolic stochastic partial differential equations.
\newblock {\em Stochastic Process. Appl.}, 123(12):4294--4336, 2013.

\bibitem{krylov-rozowskii}
N.~V. Krylov and B.~L. Rozovski\u{\i}.
\newblock Stochastic evolution equations.
\newblock In {\em Current problems in mathematics, {V}ol. 14 ({R}ussian)},
  Itogi Nauki i Tekhniki, pages 71--147, 256. Akad. Nauk SSSR, Vsesoyuz. Inst.
  Nauchn. i Tekhn. Inform., Moscow, 1979.

\bibitem{L:09}
W.~Liu.
\newblock On the stochastic {$p$}-{L}aplace equation.
\newblock {\em J. Math. Anal. Appl.}, 360(2):737--751, 2009.

\bibitem{LT:11}
W.~Liu and J.~M. T\"{o}lle.
\newblock Existence and uniqueness of invariant measures for stochastic
  evolution equations with weakly dissipative drifts.
\newblock {\em Electron. Commun. Probab.}, 16:447--457, 2011.

\bibitem{N:21}
M.~Neu{\ss}.
\newblock Well-posedness of {SVI} solutions to singular-degenerate stochastic
  porous media equations arising in self-organized criticality.
\newblock {\em Stoch. Dyn.}, 21(5):2150029, 2021.

\bibitem{N:23}
M.~Neu{\ss}.
\newblock Ergodicity for singular-degenerate stochastic porous media equations.
\newblock {\em J. Dynam. Differential Equations}, 35(2):1561--1584, 2023.

\bibitem{pardoux}
E.~Pardoux.
\newblock {\em Équations aux dérivées partielles stochastiques non
  linéaires monotones : étude de solutions fortes de type Ito}.
\newblock PhD thesis, Université Paris-Sud, 1975.

\bibitem{strong2}
G.~D. Prato, M.~Röckner, B.~L. Rozovski\u{\i}, and F.-Y. Wang.
\newblock Strong solutions of stochastic generalized porous media equations:
  Existence, uniqueness, and ergodicity.
\newblock {\em Commun. Partial Differ. Equ.}, 31(2):277--291, 2006.

\bibitem{PR}
C.~Pr\'{e}v\^{o}t and M.~R\"{o}ckner.
\newblock {\em A concise course on stochastic partial differential equations},
  volume 1905 of {\em Lecture Notes in Mathematics}.
\newblock Springer, Berlin, 2007.

\bibitem{RRW:2007}
J.~Ren, M.~R\"{o}ckner, and F.-Y. Wang.
\newblock Stochastic generalized porous media and fast diffusion equations.
\newblock {\em J. Differential Equations}, 238(1):118--152, 2007.

\bibitem{RWX:24}
M.~R\"{o}ckner, W.~Wu, and Y.~Xie.
\newblock Stochastic generalized porous media equations over {$\sigma $}-finite
  measure spaces with non-continuous diffusivity function.
\newblock {\em Potential Anal.}, 61(4):731--773, 2024.

\bibitem{Francesco}
M.~Röckner and F.~Russo.
\newblock Uniqueness for a class of stochastic fokker–planck and porous media
  equations.
\newblock {\em J. Evol. Equ}, 17(3):1049--1062, 2017.

\bibitem{SZ:21}
N.~Sapountzoglou and A.~Zimmermann.
\newblock Well-posedness of renormalized solutions for a stochastic
  {$p$}-{L}aplace equation with {$L^1$}-initial data.
\newblock {\em Discrete Contin. Dyn. Syst.}, 41(5):2341--2376, 2021.

\bibitem{SZ:22}
N.~Sapountzoglou and A.~Zimmermann.
\newblock Renormalized solutions for stochastic {$p$}-{L}aplace equations with
  {$L^1$}-initial data: the case of multiplicative noise.
\newblock {\em Discrete Contin. Dyn. Syst.}, 42(8):3979--4002, 2022.

\bibitem{SZ:23}
K.~Schmitz and A.~Zimmermann.
\newblock The stochastic {$p$}-{L}aplace equation on {$\mathbb{R}^d$}.
\newblock {\em Stoch. Anal. Appl.}, 41(5):892--917, 2023.

\bibitem{SST:23}
F.~Seib, W.~Stannat, and J.~M. T\"{o}lle.
\newblock Stability and moment estimates for the stochastic singular
  {$\Phi$}-{L}aplace equation.
\newblock {\em J. Differential Equations}, 377:663--693, 2023.

\bibitem{T:18}
J.~M. T\"{o}lle.
\newblock Estimates for nonlinear stochastic partial differential equations
  with gradient noise via {D}irichlet forms.
\newblock In {\em Stochastic partial differential equations and related
  fields}, volume 229 of {\em Springer Proc. Math. Stat.}, pages 249--262.
  Springer, Cham, 2018.

\bibitem{T:20}
J.~M. T\"{o}lle.
\newblock Stochastic evolution equations with singular drift and gradient noise
  via curvature and commutation conditions.
\newblock {\em Stochastic Process. Appl.}, 130(5):3220--3248, 2020.

\bibitem{Turra}
M.~Turra.
\newblock Existence and extinction in finite time for {S}tratonovich gradient
  noise porous media equations.
\newblock {\em Evol. Equ. Control Theory}, 8(4), 2019.

\bibitem{VASQUEZ}
J.~V\'{a}squez.
\newblock {\em Smoothing and Decay Estimates for Nonlinear Diffusion Equations:
  {E}quations of Porous Medium Type}.
\newblock Oxford University Press, 2006.

\bibitem{W:23}
J.~Wichmann.
\newblock On temporal regularity of strong solutions to stochastic
  {$p$}-{L}aplace systems.
\newblock {\em SIAM J. Math. Anal.}, 55(4):3713--3730, 2023.

\end{thebibliography}
\end{document}